\documentclass[acmsmall,screen]{acmart}

\AtBeginDocument{%
  \providecommand\BibTeX{{%
    \normalfont B\kern-0.5em{\scshape i\kern-0.25em b}\kern-0.8em\TeX}}}

\acmConference[Woodstock '18]{Woodstock '18: ACM Symposium on Neural
  Gaze Detection}{June 03--05, 2018}{Woodstock, NY}
\acmBooktitle{Woodstock '18: ACM Symposium on Neural Gaze Detection,
  June 03--05, 2018, Woodstock, NY}
\acmPrice{15.00}
\acmISBN{978-1-4503-XXXX-X/18/06}

\acmSubmissionID{TOCL-2022-0025}

\begin{document}

\title{Generalized Realizability and Intuitionistic Logic}

\author{Aleksandr Yu. Konovalov}
\affiliation{%
  \institution{Faculty of Mechanics and Mathematics, Lomonosov Moscow State University}
  \streetaddress{GSP-1, Leninskie Gory}
  \city{Moscow}
  \postcode{119991}
  \country{Russian Federation}
}
\email{alexandr.konoval@gmail.com}

\begin{abstract}
Let $V$ be a set of number-theoretical functions.
We define a notion of $V$-realizability for predicate formulas 
in such a way
that the indices of functions in $V$ are used for interpreting the implication and the universal quantifier.
In this paper we prove that
Intuitionistic Predicate Calculus is sound with respect to the semantics of $V$-realizability
if and only if
some natural conditions for $V$ hold.
\end{abstract}



\newcommand{\uu}{\mathsf u}
\newcommand{\pp}{\mathsf p}
\newcommand{\cc}{\mathsf c}
\newcommand{\II}{\mathsf I}
\newcommand{\NN}{\mathbb N}
\newcommand{\FF}{\mathbb F}
\newcommand{\TT}{\mathcal T}
\renewcommand{\ss}{\mathsf s}

\newcommand{\cb}{BF}
\newcommand{\cpp}{PV}
\newcommand{\cefp}{DV}
\newcommand{\cek}{Cm}
\newcommand{\ceu}{Cs}
\newcommand{\cec}{Cn}
\newcommand{\csmn}{SMN}
\newcommand{\cfp}{DV}

\newcommand{\sx}{{\overline x}}
\newcommand{\sy}{{\overline y}}
\newcommand{\sxp}{{\overline{x}'}}
\newcommand{\syp}{{\overline{y}'}}
\newcommand{\st}{{\overline t}}
\newcommand{\sz}{{\overline z}}
\newcommand{\su}{{\overline u}}
\newcommand{\sv}{{\overline v}}
\newcommand{\sr}{{\overline r}}
\newcommand{\sk}{{\overline k}}
\newcommand{\sm}{{\overline m}}
\newcommand{\sa}{{\overline a}}
\newcommand{\sd}{{\overline d}}
\renewcommand{\sb}{{\overline b}}

\newcommand{\zz}{\mathsf z}
\newcommand{\xx}{\mathsf x}
\newcommand{\yy}{\mathsf y}
\newcommand{\kk}{\mathsf k}

\newcommand{\rvf}{\mathrel{\mathbf{r}_f}}
\newcommand{\rvfp}{\mathrel{\mathbf{r}'_f}}
\newcommand{\rvfx}[1]{\mathrel{\mathbf{r}_{f,\,#1}}}

\newcommand{\BQC}{\ensuremath{\mathsf{BQC}}}
\newcommand{\IPC}{\ensuremath{\mathsf{IPC}}}

\maketitle

\section{Introduction}

In \cite{klini1945} S.~Kleene introduced the notion of recursive realizability.
This semantics of mathematical statements is the foundation of the constructive approach to mathematics.
Logical laws acceptable from the constructive point of view are of interest for the development of constructive mathematics.
In mathematical logic, logical laws are expressed by means of predicate formulas.
In recursive realizability semantics for predicate formulas,
formulas of the language of formal arithmetic are substituted for predicate variables.
However, it has been shown \cite{plisko_new} that the class of realizable predicate formulas for the language of arithmetic
becomes narrower under the extension of this language by means of a truth predicate.
The efforts to produce an adequate notion of realizability for predicate formulas
independent of the language in which the predicates substituted for predicate variables are formulated
has led to the notion of absolutely recursive realizable predicate formula \cite{plisko1984}.

Propositional and predicate logics of recursive realizability were investigated since the 50s of the last century.
Constructive logics based on the recursive realizability essential differ from classical and intuitionistic logics.
As an example, the predicate logic of recursive realizability is not recursively enumerable \cite{plisko1977},
while classical and intuitionistic predicate logics are conversely recursively enumerable.
Various forms of subrecursive realizability were considered: primitive recursive realizability \cite{dam1994,sal2001_1},
minimal realizability \cite{dam1995}, and corresponding predicate logics \cite{viter_dis,pak_dis,kon_mz_sprr}.
Since these logics are not arithmetical, it is of interest to generalize the constructive approach.
As noted by H.~Rogers in \cite[\S 16.5]{rodgers}, a hyperarithmetical computability is one of the natural extensions of the constructive approach.
Recently, V.E.~Plisco and A.Yu.~Konovalov studied
a notion of realizability based on arithmetical and hyperarithmetical computability.

In paper \cite{kon_gr_ACM} we introduced a notion of $V$-realizability to generalize
a lot of known realizability semantics for predicate formulas:
recursive realizability, primitive recursive realizability by Salehi \cite{sal2001_1},
general recursive realizability \cite{kon_grr_bl},
arithmetical realizability \cite{kon_ar_bl},
hyperarithmetical realizability \cite{kon_plisko_hr}, and other.
The semantics of $V$-realizability for some set of functions $V$
is a modification of recursive realizability such that functions from the set~$V$ are used instead of partial recursive functions.
Thus usual recursive realizability is a special case of $V$-realizability, when $V$ is the set of all partial recursive functions;
arithmetical realizability is a special case of $V$-realizability, when $V$ is the set of all arithmetical functions, etc.

It is of interest to consider how logics of $V$-realizability dependence on properties of $V$.
Intuitionistic Logic is sound with respect to the semantics of recursive realizability. 
But in general this is not the case for the $V$-realizability
\cite{kon_ar_bl,kon_plisko_hr}.
For example, the formula $\forall x\,((\top \to P(x)) \to P(x))$ is not arithmetical realizable \cite{kon_Lr_IPC}.
In \cite{kon_ar_bl,kon_plisko_hr, kon_grr_bl} we prove that Basic Predicate Calculus ($\BQC$) \cite{ruitenburg}
is sound with respect to some special cases of $V$-realizability.
In paper \cite{kon_gr_ACM} we generalize this result and find natural sufficient conditions for $V$ such that
$\BQC$ is sound with respect to the semantics of $V$-realizability.
The purpose of this article is to find an additional necessary and sufficient condition on $V$ such that
Intuitionistic Logic is sound with respect to the semantics of $V$-realizability.

\section{Definitions}
\subsection{$V$-computability}
In this section we introduce some notation from \cite{kon_gr_ACM}.
Denote by $\NN$ the set of all natural numbers $0, 1, 2, \ldots$
Let $\cc$ be an arbitrary bijection of $\NN^2$ to $\NN$. 
Denote by $\pp_1, \pp_2$ the $1$-ary functions of $\NN$ to $\NN$
such that, for all $a, b \in \NN$, $\pp_1(\cc(a, b)) = a$ and $\pp_2(\cc(a, b)) = b$.
We omit the brackets in expressions of the form $\pp_1(t'),\ \pp_2(t'')$ and write $\pp_1 t',\ \pp_2 t''$.
Suppose $n \ge 1$ and $1 \le i \le n$, denote by $I^i_n$ the function of $\NN^n$ to $\NN$ such that $I^i_n(a_1, \ldots, a_n) = a_i$
for all $a_1, \ldots, a_n \in \NN$.
Let $p$ be a $n$-ary partial function and $a_1, \ldots, a_n$ natural numbers;
then we write $!p(a_1, \ldots, a_n)$ if $p(a_1, \ldots, a_n)$ is defined.

We consider an arbitrary (countable) set $V$ of partial functions with arguments and values from $\NN$.
We say that $\varphi$ is a $V$-function if $\varphi\in V$.
For every $n \ge 0$, denote by $V_n$ the set of all $n$-ary $V$-functions.
Clearly, $V = \bigcup^\infty_{n = 0} V_n$.
For every $n \ge 0$, let us fix some numbering of the set $V_n$.
This means that we fix some set of indices $\II_n \subseteq \NN$ and a mapping $e \mapsto \varphi^n_e$ such that
$\varphi^n_e$ is an $n$-ary $V$-function whenever $e \in \II_n$ and
every $n$-ary $V$-function is $\varphi^n_e$ for some $e \in \II_n$.
We often write $\varphi_e$ instead of $\varphi^n_e$ if there is no confusion.

Let $Var = \{x_1, x_2, \ldots\}$ be a countable set of variables.
We say that an expression $t$ is a \textit{$V$-term} if
$t$ is a natural number or
$t \in Var$ or
$t$ has the form $\varphi(t_1, \ldots, t_n)$, 
where $\varphi \in V_n$ and $t_1, \ldots, t_n$ are $V$-terms, for some $n \ge 0$.
Any $V$-term without variables is called \textit{closed}.
Suppose $e$ is a natural number and $t$ is a closed $V$-term,
then the relation ``$e$ is the value of $t$'' is defined inductively by the length of $t$:
$e$ is the value of $t$ if $t$ is the natural number $e$;
$e$ is the value of $\varphi(t_1, \ldots, t_n)$ if there are natural numbers $e_1, \dots, e_n$ such that
$e_1, \dots, e_n$ are the values of $t_1, \ldots, t_n$,
$\varphi(e_1, \ldots, e_n)$ is defined, and $e = \varphi(e_1, \ldots, e_n)$.
We say that the value of a closed $V$-term $t$ is defined
if there is a natural number $e$ such that $e$ is the value of $t$.
It can be easily checked that if the value of closed $V$-term $t$ is defined,
then there exists a unique natural number $e$ such that $e$ is the value of $t$.
In this case we denote by $\overline{t}$ the value of $t$.
Suppose $t_1, t_2$ are closed $V$-terms,
we write $t_1 \simeq t_2$
if either (i) the values of $t_1$ and $t_2$ are not defined,
or (ii) the values of $t_1$ and $t_2$ are defined
and $\overline{t}_1 = \overline{t}_2$.
Let $k_1, \ldots, k_n$ be natural numbers, $x_1, \ldots, x_n$ distinct variables, and $t$ an $V$-term,
denote by $[k_1, \ldots, k_n/x_1, \ldots, x_n]\,t$
the result of substituting $k_1, \ldots, k_n$ for all occurrences of $x_1, \ldots, x_n$ in $t$.
Suppose $t_1,\ t_2$ are $V$-terms
and all variables in $t_1$ and $t_2$ are in a list of distinct variables $x_1, \ldots, x_n$,
we write $t_1 \simeq t_2$
if for all natural numbers $k_1, \ldots, k_n$ we have
$[k_1, \ldots, k_n/x_1, \ldots, x_n]\,t_1 \simeq [k_1, \ldots, k_n/x_1, \ldots, x_n]\,t_2.$

We say that $V$ is \textit{a basic computability model} if the following conditions hold:
\begin{itemize}
  \item[(\cb)]
  $I^i_n$, $\cc$, $\pp_1$, $\pp_2$ are $V$-functions for all $n \ge 1$,\ $1 \le i \le n$;

  \item[(\cek)]
  the composition of $V$-functions is a $V$-function and
  an index of it can be obtained by some $V$-function:
  for all natural numbers $n, m_1, \ldots, m_n$ there is an $(n+1)$-ary $V$-function $s$ such that
  $s(e, e_1, \ldots, e_n) \in \II_m$ and
\begin{equation*}
  \varphi_{s(e, e_1, \ldots, e_n)}(x_1, \ldots, x_m) \simeq
  \varphi_e(\varphi_{e_1}(x_1, \ldots, x_{m_1}), \ldots, \varphi_{e_n}(x_1, \ldots, x_{m_n}))
\end{equation*}
for all~$e \in \II_n, e_1 \in \II_{m_1}, \ldots, e_n \in \II_{m_n}$,
where $m = \max_{1 \le i \le n}{m_i}$;

\item[(\cec)]
  every constant function is a $V$-function and
  an index of it can be obtained by some $V$-function:
  there exists a $V$-function $s$ such that, for all natural numbers $k$,
  we have $s(k) \in \II_0$ and $\varphi^0_{s(k)} \simeq k$;

\item[(\ceu)]
an index of a ``conditional function'' can be obtained by some $V$-function:
for every natural number $n$ there is a $V$-function $s$
such that, for all natural numbers $d$ and $e_1, e_2 \in \II_n$, we have $s(e_1, e_2) \in \II_{n+1}$,
\begin{align*}
\varphi_{s(e_1, e_2)}(x_1, \ldots, x_n, d) \simeq \varphi_{e_1}(x_1, \ldots, x_n)\ \mbox{ if } d = 0, \\
 \varphi_{s(e_1, e_2)}(x_1, \ldots, x_n, d) \simeq \varphi_{e_2}(x_1, \ldots, x_n)\ \mbox{ if } d \not= 0.
\end{align*}
\end{itemize}

For example,
if $\cc$, $\pp_1$, $\pp_2$ are recursive (see \S 5.3 in \cite{rodgers}), then
the following sets of functions with some numbering satisfy the conditions
(\cb), (\cek), (\cec), (\ceu):
\begin{itemize}
  \item the set of all partial recursive functions;
  \item the set of all total recursive functions (see \cite{kon_grr_bl, kon_grr_il});
  \item the set of all arithmetical functions (see \cite{kon_ar_bl,kon_ar_pr});
  \item the set of all hyperarithmetical functions (see \cite{kon_plisko_hr});
  \item the set of all $L$-defined functions, where $L$ is an extension of the language of arithmetic
        (see \cite{kon_Lr_IPC,kon_Lr_IPC_IS}).
\end{itemize}

Suppose $V$ is a basic computability model;
then the following conditions hold (see \cite{kon_gr_ACM}):
\begin{itemize}
  \item[(\cpp)]
  any permutation of variables is available for the $V$-functions:
  if $p$ is a permutation of the set $\{1, \ldots, n\}$,
  then there is a $V$-function $s$ such that,
  for all $e \in \II_n$, $s(e) \in \II_n$ and
  $
  \varphi_{s(e)}(x_1, \ldots, x_n) \simeq \varphi_e(x_{p(1)}, \ldots, x_{p(n)});
  $

  \item[(\cefp)] adding of a dummy variable is available for the $V$-functions:
  for all natural numbers $n$ there exists a $V$-function $s$ such that,
    for all $e \in \II_n$, $s(e) \in \II_{n+1}$ and
    $
    \varphi_{s(e)}(x_1, \ldots, x_n, x_{n+1}) \simeq \varphi_e(x_1, \ldots, x_n);
    $
  \item[(\csmn)] an analog of the ($s-m-n$)-theorem (Theorem V \S 1.8 in \cite{rodgers}):
  for all natural numbers $m,\ n$ there exists a $V$-function $s$ such that,
for all natural numbers $k_1, \ldots, k_m$ and $e \in \II_{m+n}$, we have $s(e, k_1, \ldots, k_m) \in \II_n$ and
\begin{equation*}
  \varphi_{s(e, k_1, \ldots, k_m)}(x_1, \ldots, x_n) \simeq \varphi_e(x_1, \ldots, x_n, k_1, \ldots, k_m).
\end{equation*}

\item[(\ceu$'$)]
for every natural number $n$ there is a $V$-function $s$
such that, for all natural numbers $d$ and $e_1, e_2 \in \II_{n+1}$, we have $s(e_1, e_2) \in \II_{n+1}$,
\begin{align*}
\varphi_{s(e_1, e_2)}(x_1, \ldots, x_n, d) \simeq \varphi_{e_1}(x_1, \ldots, x_n, \pp_2 d)\ \mbox{ if } \pp_1 d = 0, \\
 \varphi_{s(e_1, e_2)}(x_1, \ldots, x_n, d) \simeq \varphi_{e_2}(x_1, \ldots, x_n, \pp_2 d)\ \mbox{ if } \pp_1 d \not= 0;
\end{align*}
\end{itemize}


We say that an $(n+1)$-ary partial function $\uu^n$ is \textit{overuniversal} for the set of all $n$-ary $V$-functions
if we have $!\uu^n(e, a_1, \ldots, a_n)$ and $\uu^n(e, a_1, \ldots, a_n) = \varphi_e(a_1, \ldots, a_n)$
for all natural numbers $e, a_1, \ldots, a_n$ such that $e \in \II_n$ and $!\varphi_e(a_1, \ldots, a_n)$.

We say that $V$ is \textit{an intuitionistic computability model} if $V$ is a basic computability model
and the following condition holds:
\begin{itemize}
  \item[(U)] there exists an overuniversal $V$-function for the set of all unary $V$-functions.
\end{itemize}

Let $V$ be an intuitionistic computability model. We shall show that for every $n \ge 1$ the following condition holds:
\begin{itemize}
        \item[($\mathrm{U^n}$)] there exists an overuniversal $V$-function for the set of all $n$-ary $V$-functions.
\end{itemize}
\begin{lemma}
    (\cb), (\cek), (U) imply ($U^n$) for every $n \ge 1$.
\end{lemma}
\begin{proof}
  The proof is by induction on n. For $n = 1$ there is nothing to prove.
  Let $n \ge 1$ and $\uu^n$ be an overuniversal $V$-function for the set of all $n$-ary $V$-functions.
  By (\cek) and (\cb), there is a $V$-function $\ss$ such that for every $e \in \II_{n+1}$ we have $s(e) \in \II_{n}$ and
  \begin{equation}\label{eq:lm1:1}
    \varphi_{\ss(e)}(x_1, \ldots, x_{n-1}, y) \simeq \varphi_e(x_1, \ldots, x_{n-1}, \pp_1 y, \pp_2 y).
  \end{equation}
  It follows from \eqref{eq:lm1:1} that
  \begin{equation}\label{eq:lm1:2}
    \varphi_{\ss(e)}(x_1, \ldots, x_{n-1}, \cc(x_n, x_{n+1})) \simeq \varphi_e(x_1, \ldots, x_{n-1}, x_n, x_{n+1}).
  \end{equation}
  By (\cek) and (\cb), there is a $V$-function $\uu^{n+1}$ such that
  \begin{equation}\label{eq:lm1:3}
    \uu^{n+1}(y, x_1, \ldots, x_{n+1}) \simeq \uu^n(\ss(y), x_1, \ldots, x_{n-1}, \cc(x_n, x_{n+1}))
  \end{equation}
  Suppose $e, a_1, \ldots, a_{n+1}$ are natural numbers such that $e \in \II_{n+1}$ and $!\varphi_e(a_1, \ldots, a_{n+1})$.
  Then $s(e) \in \II_n$ and it follows from \eqref{eq:lm1:2} that $!\varphi_{\ss(e)}(a_1, \ldots, a_{n-1}, \cc(a_n, a_{n+1}))$
  and
  \begin{equation}\label{eq:lm1:21}
    \varphi_{\ss(e)}(a_1, \ldots, a_{n-1}, \cc(a_n, a_{n+1})) = \varphi_e(a_1, \ldots, a_{n+1}).
  \end{equation}
  Since $!\varphi_{\ss(e)}(a_1, \ldots, a_{n-1}, \cc(a_n, a_{n+1}))$,
  we have
  \begin{equation}\label{eq:lm1:4}
    \varphi_{\ss(e)}(a_1, \ldots, a_{n-1}, \cc(a_n, a_{n+1})) = \uu^n(\ss(e), a_1, \ldots, a_{n-1}, \cc(a_n, a_{n+1})).
  \end{equation}
  It follows from \eqref{eq:lm1:3}, \eqref{eq:lm1:21}, \eqref{eq:lm1:4} that $!\uu^{n+1}(e, a_1, \ldots, a_{n+1})$
  and
  $\uu^{n+1}(e, a_1, \ldots, a_{n+1}) = \varphi_e(a_1, \ldots, a_{n+1})$.
  Thus $\uu^{n+1}$ is an overuniversal $V$-function for the set of all $(n+1)$-ary $V$-functions.
\end{proof}

\subsection{Intuitionistic Predicate Calculus}

We consider the language of Intuitionistic Predicate Calculus ($\IPC$) without constants and functional symbols.
\textit{The language of $\IPC$} contains
a countably infinite set of predicate symbols for each finite arity,
a countably infinite set of variables,
parentheses,
the logical constants $\bot$ (falsehood), $\top$ (truth),
the logical connectives $\land$, $\lor$, $\to$
and the quantifiers $\forall$, $\exists$.
Suppose $M \subseteq \NN$, denote by $L^M_\IPC$
the extension of the language of $\IPC$ by individual constants from the set $M$.
Thus the language of $\IPC$ is a special case of $L^M_\IPC$ for $M = \varnothing$.
We write $L_\IPC$ instead of $L^\varnothing_\IPC$.

\textit{Terms} of $L^M_\IPC$ are constants from $M$ and variables.
\textit{Atoms} of $L^M_\IPC$ are $\bot,\ \top$, and expressions of the form $P(t_1, \ldots, t_n)$,
where $P$ is an $n$-ary predicate symbol
and $t_1, \ldots, t_n$ are terms of $L^M_\IPC$.
\textit{Formulas} of $L^M_\IPC$ are built up 
according to the following grammar:
\begin{equation*}\label{grammarLBQC}
  A,\, B ::= At \mid A \land B \mid A \lor B \mid A \to B \mid \forall y\, A \mid \exists y\, A;
\end{equation*}
here
$At$ is an atom of $L^M_\IPC$
and $y$ is a variable.
We write $\forall \sx\, B$ instead of $\forall x_1, \ldots, \forall x_n\, B$ for $\sx = x_1, \ldots, x_n$.
Terms and formulas of $L^M_\IPC$ will be called \textit{$M$-terms} and \textit{$M$-formulas}, for short.
At the same time formulas of $L_\IPC$ are said to be \textit{formulas}.


Free and bound variables are defined in the usual way.
An occurrence of a variable $x$ in an $M$-formula $A$ is \textit{free}
if it is not in the scope of a quantifier $\exists x$ or $\forall x$ in $A$.
An occurrence of a variable in an $M$-formula that is not free is called \textit{bound}.
We say that a variable $x$ is a \textit{free variable} (\textit{bound variable}) of an $M$-formula $A$
if there exists a free (bound) occurrence of $x$ in $A$.
A sentence of $L^M_\IPC$ is a formula of $L^M_\IPC$ without free variables.
Sentences of $L^M_\IPC$ are called \textit{$M$-sentences}, and sentences of $L_\IPC$ simply \textit{sentences}, for short.

An $M$-term $t$ is called \textit{free} for a variable $x$ in a $M$-formula $A$
if for each variable $y$ in $t$ there is no occurrence of $x$
in the scope of a quantifier $\exists y$ or $\forall y$.
Let $t_1, \ldots, t_n$ be $M$-terms, $x_1, \ldots, x_n$ be distinct variables,
and $A$ be an $M$-formula,
denote by $[t_1, \ldots, t_n/x_1, \ldots, x_n] A$
the result of substituting $t_1, \ldots, t_n$ for all free occurrences of $x_1, \ldots, x_n$
in a formula $A'$
obtained from $A$ by renaming all bound variables
in such a way that, for each $i = 1, \ldots, n$, the $M$-term $t_i$ is free for $x_i$ in $A'$.

Suppose $A$ is an $M$-formula
and all free variables of $A$ are in $\sx$,
where $\sx$ is a list of distinct variables.
By the statement ``$A(\sx)$ is a $M$-formula''
we mean the conjunction of statements: ``$A$ is an $M$-formula'',
``$\sx$ is a list of distinct variables'',
and ``all free variables of $A$ are in $\sx$''.

If $\st = t_1, \ldots, t_n$ is a list of $M$-terms,
then put $|\st| \rightleftharpoons n$.
Let $A(\sx)$ be an $M$-formula
and $\st$ be a list of $M$-terms such that $|\st| = |\sx|$;
then by $A(\st)$ denote $[\st/ \sx] A$.


The axioms of $\IPC$ are:
\begin{itemize}
\item[A1)] $\top$;
\item[A2)] $A\to(B\to A)$;
\item[A3)] $(A\to(B\to C))\to((A\to B)\to(A\to C))$;
\item[A4)] $A\to(B\to A\land B)$;
\item[A5)] $A\land B \to A$;
\item[A6)] $A\land B \to B$;
\item[A7)] $(A\to C)\to((B\to C)\to(A\lor B\to C))$;
\item[A8)] $A\to A\lor B$;
\item[A9)] $B\to A\lor B$;
\item[A10)] $\bot \to A$;
\item[A11)] $\forall y\, A \to [t / y]\,A$;
\item[A12)] $[t / y]\,A \to \exists y\,A$;
\item[A13)] $\forall x\,(B \to A)\to (B\to \forall x\,A)$ if $x$ is not free in $B$;
\item[A14)] $\forall x\,(A\to B)\to(\exists x\,A \to B)$ if $x$ is not free in $B$.
\end{itemize}

The rules of $\IPC$ are:\smallskip
\begin{itemize}
  \item[R1)] $\frac{\displaystyle A,\ A\to B}{\displaystyle B}$;\smallskip
  \item[R2)] $\frac{\displaystyle A}{\displaystyle  \forall y\,A}$.
\end{itemize}
In the axioms and rules of $\IPC$\ $A,\ B,\ C$ are formulas, $t$ is a term,
and $y$ is a variable.

Given a sequent $A$,
we write $\IPC \vdash A$ if $A$ is derivable in $\IPC$.

\subsection{$V$-realizability}


In \cite{kon_gr_for_lang_ar,kon_gr_for_lang_ar_IS} we introduced a notion of $V$-realizability for the language of arithmetic.
Using methods of \cite{plisko1984,kon_MP,kon_MP_IS},
in \cite{kon_gr_ACM}
we defined a notion of absolute $V$-realizability in some domain $M \subseteq \NN$
for the formulas of the language of Basic Predicate Calculus ($\BQC$) \cite{ruitenburg, ruitenburg_FST}.
The language of $\BQC$ ($L_\BQC$) is not the same as the usual one for $\IPC$.
The language of $\BQC$ is differ from this language by the way of using the universal quantifier.
Namely, the quantifier $\forall$ is used only in the formulas of the form $\forall \sx\,(A \to B)$,
where $\sx$ is a finite list of variables, $A$ and $B$ being formulas.
In this section we extend the definition of the notion of absolute $V$-realizability from the formulas of $L_\BQC$
to the formulas of $L_\IPC$.

Suppose $\varnothing \not= M \subseteq \NN$,
we call any total function from $M^n$ to $2^\NN$ an \textit{$n$-ary generalized predicate} on $M$,
where $2^\NN$ is the set of all subsets of $\NN$.
A mapping $f$ is called an \textit{$M$-evaluation} if
$f(P)$ is an $n$-ary generalized predicate on $M$
whenever $P$ is an $n$-ary predicate symbol of $L_\IPC$.
We write $P^f$ instead of $f(P)$.
We say that $f$ is an \textit{evaluation}
if $f$ is an $M$-evaluation for some $\varnothing \not= M \subseteq \NN$.
We say that $M$ is \textit{domain} of an evaluation $f$ if
$f$ is an $M$-evaluation.

\begin{definition}
Let $e$ be a natural number, $M$ a nonempty subset of $\NN$, $f$ an $M$-evaluation, and $A$ an $M$-sentence.
The relation  ‘‘$e$ $V$-\textit{realizes} $A$ \textit{on} $f$'' is denoted
$e \rvf A$ and
is defined by induction on the number of logical connectives
and quantifiers in $A$:
\begin{itemize}
\item there is no $e$ such that $e \rvf \bot$;
\smallskip

\item $e \rvf \top$ for all $e$;
\smallskip

\item $e \rvf P(a_1, \ldots, a_n) \rightleftharpoons e \in P^f(a_1, \ldots, a_n)$,
where $P$ is an $n$-ary predicate symbol
and $a_1, \ldots, a_n \in M$;
\smallskip

\item $e \rvf (\Phi \land \Psi) \rightleftharpoons$ $\pp_1 e \rvf \Phi$ and $\pp_2 e \rvf \Psi$;
\smallskip

\item $e \rvf (\Phi \lor \Psi) \rightleftharpoons$
$(\pp_1 e = 0$ and $\pp_2 e \rvf \Phi)$
or $(\pp_1 e = 1$ and $\pp_2 e \rvf \Psi)$;
\smallskip

\item $e \rvf \Phi \to \Psi \rightleftharpoons$
$e \in \II_1$ and,
for every $s \in \NN$,
if $s \rvf \Phi$, then
$!\varphi_e(s)$ 
and
$\varphi_e(s)\rvf \Psi$;
\smallskip

\item $e \rvf \exists x \:\Phi(x) \rightleftharpoons \pp_1 e \in M$ and $\pp_2 e \rvf \Phi(\pp_1 e)$;
\smallskip


\item $e \rvf \forall \sx\,(\Phi(\sx) \to \Psi(\sx)) \rightleftharpoons$
$e \in \II_{n+1}$  and,
for every $s \in \NN$, for all $\sa \in M$, 
if $s \rvf \Phi(\sa)$, then
$!\varphi_e(\sa, s)$ 
and
$\varphi_e(\sa, s) \rvf \Psi(\sa)$,
where $n = |\sx| = |\sa|$ and $n \ge 1$.
\smallskip

\item $e \rvf \forall \sx\,\Psi(\sx) \rightleftharpoons$
$e \rvf \forall \sx\,(\top \to \Psi(\sx))$,
if $\Psi(\sx)$ has not the form $A \to B$ or $\forall y\,B$.
\end{itemize}
\end{definition}
A sentence $A$ is called \textit{absolutely $V$-realizable over all domains}
if there exists a natural number $e$ such that 
$e \rvf A$ whenever $f$ is an evaluation. 
We say that a sentence $A$ is \textit{weak $V$-realizable} in a domain $\varnothing \not= M \subseteq \NN$
if, for every $M$-evaluatin $f$, there is a natural number $e$ such that $e \rvf A$.

\section{Main result}

Let $e$ be a natural number, $M$ a nonempty subset of $\NN$, $f$ a $M$-evaluation, and  $A(y)$ a $M$-formula;
then we write $e \rvfp \forall y\,A(y)$ if
$e \in \II_1$ and, for every $a \in M$,
$!\varphi_e(a)$ and $\varphi_e(a) \rvf A(a)$.

\begin{lemma}\label{l:translation_to_simple_forall}
    Let $V$ be an intuitionistic computability model.
    Suppose $A(y, \sz)$ is a formula, where $\sz = z_1, \ldots, z_m$;
    then there are an unary $V$-functions $g_A, h_A$ such that,
    for every natural number $e$,
    for every evaluation $f$,
    for all $\sb = b_1, \ldots, b_m \in M$ (here $M$ is the domain of $f$), we have:
    \begin{itemize}
      \item[a)] if $e \rvf \forall y\,A(y, \sb)$, then $!g_A(e)$ and $g_A(e) \rvfp \forall y\,A(y, \sb)$;
      \item[b)] if $e \rvfp \forall y\,A(y, \sb)$, then $!h_A(e)$ and $h_A(e) \rvf \forall y\,A(y, \sb)$.
    \end{itemize}
\end{lemma}
\begin{proof}\hspace{0pt}
\begin{itemize}
\item[1)] Let $A(y, \sz)$ have the form $\forall \sx\,(B(y, \sx, \sz) \to C(y, \sx, \sz))$,
  where $\sx = x_1, \ldots, x_n$ and $n \ge 0$.
  In the case $n=0$, we assume that $A(y, \sz)$ has the form $B(y, \sz) \to C(y, \sz)$.
  By (\csmn) and (\cpp), there is a $V$-function $\kk$ such that,
  for all natural number $a$ and $e \in \II_{n+2}$,
  we have $\kk(e, a) \in \II_{n+1}$ and
  \begin{equation}\label{f:tsf:s1}
    \varphi_{\kk(e, a)}(\sx, w) \simeq \varphi_e(a, \sx, w).
  \end{equation}
  It follows from (\csmn) and (\cpp) that there exists a $V$-function $g_A$
  such that, for every natural number $e$, we have $g_A(e) \in \II_1$ and
  \begin{equation}\label{f:tsf:s2}
    \varphi_{g_A(e)}(y) \simeq \kk(e, y).
  \end{equation}
  By (\cek), there is a $V$-function $h_A$ such that,
  for every $e \in \II_1$, we have $h_A(e) \in \II_{n+2}$ and
  \begin{equation}\label{f:tsf:s4}
    \varphi_{h_A(e)}(y, \sx, w) \simeq \uu^{n+1}(\varphi_e(y), \sx, w),
  \end{equation}
  where $\uu^{n+1}$ is an overuniversal $V$-function for the set of all $(n+1)$-ary $V$-functions.

  Let $\varnothing \not= M \subseteq \NN$ and $f$ be an $M$-evaluation, and $\sb = b_1, \ldots, b_m \in M$.

  \begin{itemize}
    \item[a)]
    Suppose $e \rvf \forall y\,A(y, \sb)$ for some natural number $e$,
      that is $$e \rvf \forall y,\sx\,(B(y, \sx, \sb) \to C(y, \sx, \sb)).$$
      Hence, for all $a, \sd = a, d_1, \ldots, d_n \in M$, for every natural number $s$, we have
      $!\varphi_e(a, \sd, s)$ and $\varphi_e(a, \sd, s) \rvf C(a, \sd, \sb)$
      whenever $s \rvf B(a, \sd, \sb)$.
      Therefore it follows from \eqref{f:tsf:s1} that,
      for every $a \in M$, for all $\sd = d_1, \ldots, d_n \in M$, for every natural number $s$, we have
      $!\varphi_{\kk(e, a)}(\sd, s)$ and $\varphi_{\kk(e, a)}(\sd, s) \rvf C(a, \sd, \sb)$
      whenever $s \rvf B(a, \sd, \sb)$.
      Thus, for every $a \in M$, we have
      \begin{equation}\label{f:tsf:d1}
        \kk(e, a) \in \forall \sx\,(B(a, \sx, \sb) \to C(a, \sx, \sb)).
      \end{equation}
      It follows from \eqref{f:tsf:s2}, \eqref{f:tsf:d1} that
      $\varphi_{g_A(e)}(a) \rvf A(a, \sb)$ for every $a \in M$.
      Thus $g_A(e) \rvfp \forall y\,A(y, \sb)$.
    \item[b)]
    Suppose $e \rvfp \forall y\,A(y, \sb)$ for some natural number $e$,
    that is $$e \rvfp \forall y\, \forall \sx\,(B(y, \sx, \sb) \to C(y, \sx, \sb)).$$
    Hence, for every $a \in M$, we have $!\varphi_e(a)$ and
    $\varphi_e(a) \rvf \forall \sx\,(B(a, \sx, \sb) \to C(a, \sx, \sb)).$
    Therefore, for every $a \in M$, for all $\sd = d_1, \ldots, d_n \in M$, for every natural number $s$, we have
    $!\varphi_{\varphi_e(a)}(\sd, s)$ and $\varphi_{\varphi_e(a)}(\sd, s) \rvf C(a, \sd, \sb)$
    whenever $s \rvf B(a, \sd, \sb)$.
    By definition, we get $$\uu^{n+1}(\varphi_e(a), \sd, s) = \varphi_{\varphi_e(a)}(\sd, s)$$
    if $!\varphi_{\varphi_e(a)}(\sd, s)$.
    Hence it follows from \eqref{f:tsf:s4} from that,
    for all $a, \sd = a, d_1, \ldots, d_n \in M$, for every natural number $s$, we have
    $!\varphi_{h_A(e)}(a, \sd, s)$ and $\varphi_{h_A(e)}(a, \sd, s) \rvf C(a, \sd, \sb)$
    whenever $s \rvf B(a, \sd, \sb)$.
    Thus $h_A(e) \rvf \forall y, \sx\,(B(y, \sx, \sb) \to C(y, \sx, \sb))$,
    that is $h_A(e) \rvf \forall y\,A(y, \sb)$.
  \end{itemize}

  \item[2)] Let $A(y, \sz)$ have not the form $B \to C$ or $\forall x\,B$.
  It follows from (\cek) and (\cec) that there exists a $V$-function $g_A$
  such that, for every $e \in \II_2$, we have $g_A(e) \in \II_1$ and
  \begin{equation}\label{f:tsf:g2}
    \varphi_{g_A(e)}(y) \simeq \varphi_e(y, 0).
  \end{equation}

  By (\cfp), there is a $V$-function $h_A$ such that,
  for every $e \in \II_1$, we have $h_A(e) \in \II_2$ and
  \begin{equation}\label{f:tsf:h2}
    \varphi_{h_A(e)}(y, w) \simeq \varphi_e(y).
  \end{equation}

  Let $\varnothing \not= M \subseteq \NN$ and $f$ be an $M$-evaluation, and $\sb = b_1, \ldots, b_m \in M$.
  \begin{itemize}
    \item[a)] Suppose $e \rvf \forall y\,A(y, \sb)$ for some natural number $e$.
    Therefore $e \rvf \forall y\,(\top \to A(y, \sb))$.
    Hence, for every $a \in M$, for every natural number $s$, we have $\varphi_e(a, s) \rvf A(a, \sb)$.
    Using \eqref{f:tsf:g2}, we get $g_A(e) \rvfp \forall y\,A(y, \sb)$.
    \item[b)] Suppose $e \rvfp \forall y\,A(y, \sb)$ for some natural number $e$.
    Hence $e \in \II_1$ and, for every $a \in M$, $!\varphi_e(a)$ and $\varphi_e(a) \rvf A(a, \sb)$.
    Using \eqref{f:tsf:h2}, it is easily shown that $h_A(e) \rvf \forall y\,A(y, \sb)$.
  \end{itemize}
\end{itemize}
\end{proof}

\begin{proposition}\label{pr:IPC_cor_Vreal_tex}
    Let $V$ be an intuitionistic computability model.
    Suppose $\Phi(\sz)$ is a formula 
    such that $\IPC \vdash \Phi(\sz)$, where $\sz = z_1, \ldots, z_m$;
    then there is $m$-ary $V$-function $\psi_\Phi$ such that,
    for every evaluation~$f$,
    for all $\sd = d_1, \ldots, d_m \in M$ (here $M$ is the domain of $f$),
    we have $!\psi_\Phi(\sd)$ and $\psi_\Phi(\sd) \rvf \Phi(\sd)$.
\end{proposition}
\begin{proof}
  By induction on derivations of $\Phi(\sz)$.
  Suppose $\Phi(\sz)$ is an axiom of $\IPC$.
  \begin{itemize}
  \item[A1)] Let $\Phi(\sz)$ be $\top.$
  By (\cec) and (\cfp), there is a $V$-function $\psi_\Phi$ such that $\psi_\Phi(\sz) \simeq 0.$
  Let $\varnothing \not= M \subseteq \NN$ and $f$ be an $M$-evaluation.
  It is obvious that $\psi_\Phi(\sd) \rvf \Phi(\sd)$ for all $\sd = d_1, \ldots, d_m \in M$.

  \item[A2)] Let $\Phi(\sz)$ be $A(\sz) \to (B(\sz) \to A(\sz))$.
  By (\cec) and (\cfp), there is a $V$-function $\ss$ such that, for every natural number $a$,
  we have $\ss(a) \in \II_1$ and $$\varphi_{\ss(a)}(y) \simeq a.$$
  It follows from (\cfp) that there exists a $V$-function $\psi_\Phi$ such that $$\psi_\Phi(\sz) \simeq e,$$
  where $e$ is an index of $\ss$, that is $\varphi_e(x) \simeq \ss(x)$.

  Let $\varnothing \not= M \subseteq \NN$ and $f$ be an $M$-evaluation, and $\sd = d_1, \ldots, d_m \in M$.
  Suppose $a \rvf A(\sd)$.
  For every natural number $b$ such that $b \rvf B(\sd)$, we have $\varphi_{\ss(a)}(b) = a \rvf A(\sd)$.
  Thus $\ss(a) \rvf B(\sd) \to A(\sd)$ for every $a$ such that $a \rvf A(\sd)$.
  Therefore $e \rvf A(\sd) \to (B(\sd) \to A(\sd))$.
  Hence $\psi_\Phi(\sd) \rvf \Phi(\sd)$ for all $\sd = d_1, \ldots, d_m \in M$.

  \item[A3)] Let $\Phi(\sz)$ be $(A(\sz)\to(B(\sz)\to C(\sz)))\to((A(\sz)\to B(\sz))\to(A(\sz)\to C(\sz)))$.
  By (\cek), there is a $V$-function $\ss$ such that, for all $p, q \in \II_1$, we have
  $\ss(p, q) \in \II_1$ and $$\varphi_{\ss(p, q)}(x) \simeq \uu(\varphi_p(x), \varphi_q(x)),$$
  where $\uu$ is a overuniversal $V$-function for all unary $V$-functions.
  By (\csmn) and (\cpp), there exists a $V$-function $\ss'$ such that, for every natural number $p$,
  we have $\ss'(p) \in \II_1$ and $$\varphi_{\ss'(p)}(y) \simeq \ss(p, y).$$
  It follows from (\cfp) that there exists a $V$-function $\psi_\Phi$ such that $\psi_\Phi(\sz) \simeq e$,
  where $e$ is an index of $\ss'$, that is $\varphi_e(x) \simeq \ss'(x)$.

  Let $\varnothing \not= M \subseteq \NN$ and $f$ be an $M$-evaluation, and $\sd = d_1, \ldots, d_m \in M$.
  Suppose 
  $$p \rvf A(\sd)\to(B(\sd)\to C(\sd)),\quad q \rvf A(\sd)\to B(\sd),\quad a \rvf A(\sd).$$
  It follows that $\varphi_p(a) \rvf B(\sd)\to C(\sd)$ and $\varphi_q(a) \rvf B(\sd)$.
  Thus $$\uu(\varphi_p(a), \varphi_q(a)) = \varphi_{\varphi_p(a)}(\varphi_q(a)) \rvf C(\sd)$$
  for every $a$ such that $a \rvf A(\sd)$.
  Hence $\varphi_{\ss'(p)}(q) = \ss(p, q) \rvf A(\sd) \to C(\sd)$
  for every $q$ such that $q \rvf A(\sd)\to B(\sd)$.
  Thus $$\varphi_e(p) = \ss'(p) \rvf (A(\sd)\to B(\sd)) \to (A(\sd) \to C(\sd))$$
  for every $p$ such that $p \rvf A(\sd)\to(B(\sd)\to C(\sd))$.
  Hence $\psi_\Phi(\sd) \rvf \Phi(\sd)$ for all $\sd = d_1, \ldots, d_m \in M$.

  \item[A4)] Let $\Phi(\sz)$ be $A(\sz)\to(B(\sz)\to A(\sz) \land B(\sz))$.
  It follows from (\cb), (\csmn), (\cpp) that there exists a $V$-function $\ss$ such that, for every natural number $a$,
  we have $\ss(a) \in \II_1$ and $$\varphi_{\ss(a)}(y) \simeq \cc(a, y).$$
  By (\cfp), there exists a $V$-function $\psi_\Phi$ such that $$\psi_\Phi(\sz) \simeq e,$$
  where $e$ is an index of $\ss$, that is $\varphi_e(x) \simeq \ss(x)$.

  Let $\varnothing \not= M \subseteq \NN$ and $f$ be an $M$-evaluation, and $\sd = d_1, \ldots, d_m \in M$.
  Suppose $a \rvf A(\sd)$.
  Hence $$\varphi_{\ss(a)}(b) = \cc(a, b) \rvf A(\sd) \land B(\sd)$$ for every $b$ such that $b \rvf B(\sd)$.
  Thus $$\ss(a) \rvf B(\sd) \to A(\sd) \land B(\sd)$$ for every $a$ such that $a \rvf A(\sd)$.
  It follows that $e \rvf A(\sd) \to (B(\sd) \to A(\sd) \land B(\sd))$.
  Hence $\psi_\Phi(\sd) \rvf \Phi(\sd)$ for all $\sd = d_1, \ldots, d_m \in M$.

  \item[A5)] Let $\Phi(\sz)$ be $A(\sz)\land B(\sz) \to A(\sz)$.
  By (\cfp), there exists a $V$-function $\psi_\Phi$ such that $\psi_\Phi(\sz) \simeq e$,
  where $e$ is an index of $\pp_1$, that is $\varphi_e(x) \simeq \pp_1 x$.
  It can easily be checked that, for every an evaluation $f$ with domain $M$,
  we have $\psi_\Phi(\sd) \rvf \Phi(\sd)$ for all $\sd = d_1, \ldots, d_m \in M$.

  \item[A6)] Let $\Phi(\sz)$ be $A(\sz)\land B(\sz) \to B(\sz)$.
  By (\cfp), there exists a $V$-function $\psi_\Phi$ such that $\psi_\Phi(\sz) \simeq e$,
  where $e$ is an index of $\pp_2$, that is $\varphi_e(x) \simeq \pp_2 x$.
  It can easily be checked that, for every an evaluation $f$ with domain $M$,
  we have $\psi_\Phi(\sd) \rvf \Phi(\sd)$ for all $\sd = d_1, \ldots, d_m \in M$.

  \item[A7)] Let $\Phi(\sz)$ be $(A(\sz)\to C(\sz))\to((B(\sz)\to C(\sz))\to(A(\sz)\lor B(\sz)\to C(\sz)))$.
  It follows from (\ceu$'$), (\cb) that there is a $V$-function $\ss$ such that,
  for all $p, q \in \II_1$, we have $\ss(p, q) \in \II_1$ and
  \begin{align*}
    \varphi_{\ss(p, q)}(k) \simeq \varphi_p(\pp_2 k) \text{ if } \pp_1 k = 0; \\
    \varphi_{\ss(p, q)}(k) \simeq \varphi_q(\pp_2 k) \text{ if } \pp_1 k \not= 0.
  \end{align*}
  By (\csmn) and (\cpp), there exists a $V$-function $\ss'$ such that, for every natural number $p$,
  we have $\ss'(p) \in \II_1$ and $$\varphi_{\ss'(p)}(y) \simeq \ss(p, y).$$
  It follows from (\cfp) that there exists a $V$-function $\psi_\Phi$ such that $\psi_\Phi(\sz) \simeq e$,
  where $e$ is an index of $\ss'$, that is $\varphi_e(x) \simeq \ss'(x)$.

  Let $\varnothing \not= M \subseteq \NN$ and $f$ be an $M$-evaluation, and $\sd = d_1, \ldots, d_m \in M$.
  Suppose
  \begin{equation*}
    p \rvf A(\sd)\to C(\sd),\quad
    q \rvf B(\sd)\to C(\sd),\quad
    k \rvf A(\sd) \lor B(\sd).
  \end{equation*}
  It follows that $\varphi_p(\pp_2 k) \rvf C(\sd)$ if $\pp_1 k = 0$,
  and $\varphi_q(\pp_2 k) \rvf C(\sd)$ if $\pp_1 k \not= 0$.
  Thus $$\varphi_{\ss(p, q)}(k) \rvf C(\sd)$$ for every $k$ such that $k \rvf A(\sd) \lor B(\sd)$.
  Hence $$\varphi_{\ss'(p)}(q) = \ss(p, q) \rvf A(\sd) \lor B(\sd) \to C(\sd)$$ for every $q$ such that $q \rvf B(\sd)\to C(\sd)$.
  It follows that $$\varphi_e(p) = \ss'(p) \rvf (B(\sd)\to C(\sd)) \to (A(\sd) \lor B(\sd) \to C(\sd))$$
  for every $p$ such that $p \rvf A(\sd)\to C(\sd)$.
  Hence $\psi_\Phi(\sd) \rvf \Phi(\sd)$ for all $\sd = d_1, \ldots, d_m \in M$.

  \item[A8)] Let $\Phi(\sz)$ be $A(\sz)\to A(\sz)\lor B(\sz)$.
  It follows from (\cb), (\cek), (\cec) there is a natural number $e$ such that
  $\varphi_e(x) \simeq \cc(0, x)$.
  By (\cfp), there exists a $V$-function $\psi_\Phi$ such that $\psi_\Phi(\sz) \simeq e$.
  It can easily be checked that, for every an evaluation $f$ with domain $M$,
  we have $\psi_\Phi(\sd) \rvf \Phi(\sd)$ for all $\sd = d_1, \ldots, d_m \in M$.

  \item[A9)] Let $\Phi(\sz)$ be $B(\sz)\to A(\sz)\lor B(\sz)$.
  It follows from (\cb), (\cek), (\cec) there is a natural number $e$ such that
  $\varphi_e(x) \simeq \cc(1, x)$.
  By (\cfp), there exists a $V$-function $\psi_\Phi$ such that $\psi_\Phi(\sz) \simeq e$.
  It can easily be checked that, for every an evaluation $f$ with domain $M$,
  we have $\psi_\Phi(\sd) \rvf \Phi(\sd)$ for all $\sd = d_1, \ldots, d_m \in M$.

  \item[A10)] Let $\Phi(\sz)$ be $\bot \to A(\sz)$.
  It follows from (\cb) there is a natural number $e$ such that $\varphi_e(x) \simeq x$.
  By (\cfp), there exists a $V$-function $\psi_\Phi$ such that $\psi_\Phi(\sz) \simeq e$.
  It can easily be checked that, for every an evaluation $f$ with domain $M$,
  we have $\psi_\Phi(\sd) \rvf \Phi(\sd)$ for all $\sd = d_1, \ldots, d_m \in M$.

  \item[A11)] Let $\Phi(\sz)$ be $\forall y\, A(y, z_2, \ldots, z_m) \to A(z_1, z_2, \ldots, z_m)$.
  It follows from (\cek), (\csmn), (\cpp), (\cfp) that there exists a $V$-function $\psi_\Phi$ such that,
  for all natural number $d_1, \ldots, d_m$,
  we have $\psi_\Phi(d_1, \ldots, d_m) \in \II_1$ and
  \begin{equation}\label{eq:a11:defpsi}
    \varphi_{\psi_\Phi(d_1, \ldots, d_m)}(x) \simeq \uu(g_A(x), d_1),
  \end{equation}
  where $\uu$ is a overuniversal $V$-function for all unary $V$-functions and
  $g_A$ is the function from Lemma~\ref{l:translation_to_simple_forall}.

  Let $\varnothing \not= M \subseteq \NN$ and $f$ be an $M$-evaluation, and $\sd = d_1, \ldots, d_m \in M$.
  Suppose
  \begin{equation}\label{eq:a11:s0}
    p \rvf \forall y\, A(y, d_2, \ldots, d_m).
  \end{equation}
  Using Lemma~\ref{l:translation_to_simple_forall}, we get $!g_A(p)$ and
  \begin{equation}\label{eq:a11:s}
    g_A(p) \rvfp \forall y\, A(y, d_2, \ldots, d_m).
  \end{equation}
  It follows from \eqref{eq:a11:s} that $!\varphi_{g_A(p)}(d_1)$ and
  \begin{equation}\label{eq:a11:s2}
    \varphi_{g_A(p)}(d_1) \rvf A(d_1, d_2, \ldots, d_m).
  \end{equation}
  Using \eqref{eq:a11:defpsi} and \eqref{eq:a11:s2}, we get
  \begin{equation}\label{eq:a11:s3}
    \varphi_{\psi_\Phi(d_1, \ldots, d_m)}(p) \rvf A(d_1, d_2, \ldots, d_m).
  \end{equation}
  Thus for every natural number $p$ we have \eqref{eq:a11:s3} whenever \eqref{eq:a11:s0}.
  Hence $\psi_\Phi(\sd) \in \Phi(\sd)$.

  \item[A12)] Let $\Phi(\sz)$ be $A(z_1, z_2, \ldots, z_m) \to \exists y\,A(y, z_2, \ldots, z_m)$.
  It follows from (\cek), (\csmn), (\cpp), (\cfp) that there exists a $V$-function $\psi_\Phi$ such that,
  for all natural number $d_1, \ldots, d_m$,
  we have $\psi_\Phi(d_1, \ldots, d_m) \in \II_1$ and
  \begin{equation}\label{eq:a12:defpsi}
    \varphi_{\psi_\Phi(d_1, \ldots, d_m)}(x) \simeq \cc(d_1, x).
  \end{equation}
  Let $\varnothing \not= M \subseteq \NN$ and $f$ be an $M$-evaluation, and $\sd = d_1, \ldots, d_m \in M$.
  Suppose
  \begin{equation}\label{eq:a12:s0}
    p \rvf A(d_1, d_2, \ldots, d_m).
  \end{equation}
  Therefore
  \begin{equation}\label{eq:a12:s00}
    \cc(d_1, p) \rvf \exists y\,A(y, d_2, \ldots, d_m).
  \end{equation}
  It follows from \eqref{eq:a12:defpsi} and \eqref{eq:a12:s00} that
  \begin{equation}\label{eq:a12:s1}
    \varphi_{\psi_\Phi(d_1, \ldots, d_m)}(p) \rvf \exists y\,A(y, d_2, \ldots, d_m).
  \end{equation}
  Thus for every natural number $p$ we have \eqref{eq:a12:s1} whenever \eqref{eq:a12:s0}.
  Hence $\psi_\Phi(\sd) \in \Phi(\sd)$.

  \item[A13)] Let $\Phi(\sz)$ be $\forall y\,(B(\sz) \to A(y, \sz))\to (B(\sz)\to \forall y\,A(y, \sz))$.
  By (\csmn), there is a $V$-function $\ss$ such that, for every $p \in \II_2$, for every natural number $b$,
  we have $\ss(p, b) \in \II_1$ and
  \begin{equation}\label{eq:a13:defs}
    \varphi_{\ss(p, b)}(y) \simeq \varphi_p(y, b).
  \end{equation}
  By (\cek), (\csmn), (\cpp) there exists a $V$-function $\ss'$ such that, for every natural number $p$,
  we have $\ss'(p) \in \II_1$ and
  \begin{equation}\label{eq:a13:defsp}
    \varphi_{\ss'(p)}(w) \simeq h_A(\ss(p, w)),
  \end{equation}
  where $h_A$ is the function from Lemma~\ref{l:translation_to_simple_forall}.
  It follows from (\cfp) that there exists a $V$-function $\psi_\Phi$ such that $\psi_\Phi(\sz) \simeq e$,
  where $e$ is an index of $\ss'$, that is $\varphi_e(x) \simeq \ss'(x)$.

  Let $\varnothing \not= M \subseteq \NN$ and $f$ be an $M$-evaluation, and $\sd = d_1, \ldots, d_m \in M$.
  Suppose
  \begin{equation}\label{eq:a13:s0}
    p \rvf \forall y\,(B(\sd) \to A(y, \sd)).
  \end{equation}
  Let
  \begin{equation}\label{eq:a13:s00}
    b \rvf B(\sd).
  \end{equation}
  Suppose $k \in M$.
  Using \eqref{eq:a13:s0} and \eqref{eq:a13:s00}, we get
  \begin{equation}\label{eq:a13:s1}
    \varphi_p(k, b) \rvf A(k, \sd).
  \end{equation}
  It follows from \eqref{eq:a13:defs}, \eqref{eq:a13:s1} that
  \begin{equation}\label{eq:a13:s2}
     \varphi_{\ss(p, b)}(k) \rvf A(k, \sd).
  \end{equation}
  Thus for every $k \in M$ we have \eqref{eq:a13:s2}. Hence
  \begin{equation}\label{eq:a13:s3}
     \ss(p, b)\rvfp \forall y\,A(y, \sd).
  \end{equation}
  Using \eqref{eq:a13:s3} and Lemma~\ref{l:translation_to_simple_forall}, we get
  \begin{equation}\label{eq:a13:s4}
     h_A(\ss(p, b)) \rvf \forall y\,A(y, \sd).
  \end{equation}
  It follows from \eqref{eq:a13:defsp}, \eqref{eq:a13:s4} that
  \begin{equation}\label{eq:a13:s5}
     \varphi_{\ss'(p)}(b) \rvf \forall y\,A(y, \sd).
  \end{equation}
  Thus for every natural number $b$ we have \eqref{eq:a13:s5} whenever \eqref{eq:a13:s00}.
  Hence
  \begin{equation}\label{eq:a13:s6}
     \varphi_e(p) = \ss'(p) \rvf B(\sd) \to \forall y\,A(y, \sd).
  \end{equation}
  Thus for every natural number $p$ we have \eqref{eq:a13:s6} whenever \eqref{eq:a13:s0}.
  Therefor $\psi_\Phi(\sd) = e \rvf \Phi(\sd)$.

  \item[A14)] Let $\Phi(\sz)$ be $\forall y\,(A(y, \sz)\to B(\sz))\to(\exists y\,A(y, \sz) \to B(\sz))$.
  By (\cek), (\cb) and (\csmn), there is a $V$-function $\ss$ such that, for every $p \in \II_2$,
  we have $\ss(p) \in \II_1$ and
  \begin{equation}\label{eq:a14:defs}
    \varphi_{\ss(p)}(x) \simeq \varphi_p(\pp_1 x, \pp_2 x).
  \end{equation}
  It follows from (\cfp) that there exists a $V$-function $\psi_\Phi$ such that $\psi_\Phi(\sz) \simeq e,$
  where $e$ is an index of $\ss$, that is $\varphi_e(x) \simeq \ss(x)$.

  Let $\varnothing \not= M \subseteq \NN$ and $f$ be an $M$-evaluation, and $\sd = d_1, \ldots, d_m \in M$.
  Suppose
  \begin{equation}\label{eq:a14:s0}
    p \rvf \forall y\,(A(y, \sd)\to B(\sd)).
  \end{equation}
  Let
  \begin{equation}\label{eq:a14:s00}
    q \rvf \exists y\,A(y, \sd).
  \end{equation}
  Using \eqref{eq:a14:s00}, we get
  \begin{equation}\label{eq:a14:s000}
    \pp_2 q \rvf A(\pp_1 q, \sd).
  \end{equation}
  It follows from \eqref{eq:a14:s0}, \eqref{eq:a14:s000} that
  \begin{equation}\label{eq:a14:s1}
    \varphi_p(\pp_1 q, \pp_2 q) \rvf B(\sd).
  \end{equation}
  Using \eqref{eq:a14:defs} and \eqref{eq:a14:s1}, we get
  \begin{equation}\label{eq:a14:s2}
    \varphi_{\ss(p)}(q) \rvf B(\sd).
  \end{equation}
  Thus for every natural number $q$ we have \eqref{eq:a14:s2} whenever \eqref{eq:a14:s00}.
  Hence
  \begin{equation}\label{eq:a14:s3}
    \varphi_e(p) = \ss(p) \rvf \exists y\,A(y, \sd) \to B(\sd).
  \end{equation}
  Thus for every natural number $p$ we have \eqref{eq:a14:s3} whenever \eqref{eq:a14:s0}.
  Therefor $\psi_\Phi(\sd) = e \rvf \Phi(\sd)$.
  \end{itemize}
  Suppose $\Phi(\sz)$ is obtained by a rule of $\IPC$.
  \begin{itemize}
    \item[R1)] Let $\Phi(\sz)$ be obtained by
    $\frac{\displaystyle A(\sz, \sy),\ A(\sz, \sy)\to \Phi(\sz)}{\displaystyle \Phi(\sz)}$.
    By the induction hypothesis, there exist $V$-functions $\psi_A$ and $\psi_{A \to \Phi}$ such that,
    for every an evaluation $f$ with domain $M$, we have
    \begin{eqnarray*}
      &\psi_A(\sd, \sa) \rvf A(\sd, \sa), \\
      &\psi_{A \to \Phi}(\sd, \sa) \rvf A(\sd, \sa) \to \Phi(\sd)
    \end{eqnarray*}
    for all $\sd, \sa \in M$.
    It follows from (\cek), (\cb) that there is a $V$-function $\psi_\Phi$ such that
    \begin{equation}\label{eq:r1:defpsi}
      \psi_\Phi(\sz) \simeq \uu(\psi_{A \to \Phi}(\sz, z_1, \ldots, z_1),\ \psi_A(\sz, z_1, \ldots, z_1)),
    \end{equation}
    where $\uu$ is a overuniversal $V$-function for the set of all unary $V$-functions.

    Let $\varnothing \not= M \subseteq \NN$ and $f$ be an $M$-evaluation, and $\sd = d_1, \ldots, d_m \in M$.
    Then
    \begin{eqnarray*}
      &\psi_A(\sd, d_1, \dots, d_1) \rvf A(\sd, d_1, \dots, d_1), \\
      &\psi_{A \to \Phi}(\sd, d_1, \dots, d_1) \rvf A(\sd, d_1, \dots, d_1) \to \Phi(\sd).
    \end{eqnarray*}
    Hence
    \begin{equation}\label{eq:r1:s0}
      \varphi_{\psi_{A \to \Phi}(\sd, d_1, \dots, d_1)}(\psi_A(\sd, d_1, \dots, d_1)) \rvf \Phi(\sd).
    \end{equation}
    It follows from \eqref{eq:r1:defpsi}, \eqref{eq:r1:s0} that
    \begin{equation*}
      \psi_\Phi(\sd) = \uu(\psi_{A \to \Phi}(\sd, d_1, \dots, d_1),\ \psi_A(\sd, d_1, \dots, d_1))
      \rvf \Phi(\sd).
    \end{equation*}

    \item[R2)] Let $\Phi(\sz)$ be obtained by $\frac{\displaystyle A(y, \sz)}{\displaystyle  \forall y\,A(y, \sz)}$.
    By the induction hypothesis, there exists a $V$-function $\psi_A$ such that,
    for every an evaluation $f$ with domain $M$, for all $a, \sd \in M$, we have
    \begin{equation*}
      \psi_A(a, \sd) \rvf A(a, \sd).
    \end{equation*}
    By (\csmn), there exists an $m$-ary $V$-function $\ss$ such that,
    for all natural numbers $\sd = d_1, \ldots, d_m$, we have $\ss(\sd) \in \II_1$ and
    \begin{equation*}
      \varphi_{\ss(\sd)}(y) \simeq \psi_A(y, \sd).
    \end{equation*}
    It follows from (\cek) that there is a $V$-function $\psi_\Phi$ such that
    \begin{equation*}
      \psi_\Phi(\sz) \simeq h_A(\ss(\sz)),
    \end{equation*}
    where $h_A$ is the function from Lemma~\ref{l:translation_to_simple_forall}.

    Let $\varnothing \not= M \subseteq \NN$ and $f$ be an $M$-evaluation, and $\sd = d_1, \ldots, d_m \in M$.
    For every $a \in M$ we have
    \begin{equation*}
      \varphi_{\ss(\sd)}(a) = \psi_A(a, \sd) \rvf A(a, \sd).
    \end{equation*}
    Hence
    \begin{equation*}
      \ss(\sd) \rvfp \forall y\,A(y, \sd).
    \end{equation*}
    Using Lemma~\ref{l:translation_to_simple_forall}, we get
    \begin{equation*}
      \psi_\Phi(\sd) = h_A(\ss(\sd)) \rvf \forall y\,A(y, \sd).
    \end{equation*}
    Thus we have $\psi_\Phi(\sd) \rvf \Phi(\sd)$ for all $\sd = d_1, \ldots, d_m \in M$.

  \end{itemize}

\end{proof}

Proposition~\ref{pr:IPC_cor_Vreal_tex} imply the next statement.
\begin{proposition}\label{pr:IPC_is_sound}
    Let $V$ be an intuitionistic computability model.
    Suppose $\Phi$ is a sentence and $\IPC \vdash \Phi$;
    then $\Phi$ is absolutely $V$-realizable over all domains.
\end{proposition}

Our main result is the following.
\begin{theorem}\label{t_main}
  Let $V$ be an basic computability model;
  then the following conditions are equivalent:
  \begin{itemize}
    \item[(i)] Intuitionistic Predicate Calculus is sound with respect to the semantics of absolute $V$-realizability over all domains;
    \item[(ii)] Intuitionistic Predicate Calculus is sound with respect to the semantics of weak $V$-realizability in the domain $\NN$;
    \item[(iii)] the formula
                \begin{equation}\label{eq:contr_exam}
                  \forall x\,(Q(x) \to \forall y\,(R(x, y) \to \exists z\,P(x, y, z))) \to
                 \forall y\,\forall x\,(Q(x) \land R(x, y) \to \exists z\,P(x, y, z))
                \end{equation}
               is weak $V$-realizable in the domain $\NN$;
    \item[(iv)] 
                there is an overuniversal $V$-function for the set of all unary $V$-functions.
  \end{itemize}
\end{theorem}
\begin{proof}
  $\mathrm{(i)} \Rightarrow \mathrm{(ii)}$ is trivial.
  $\mathrm{(ii) \Rightarrow (iii)}$ is true,
  because \eqref{eq:contr_exam} is deducible in $\IPC$.

  Suppose $\mathrm{(iii)}$ is true.
  Let $f$ be a $\NN$-evaluation such that, for all natural numbers $e, a, b, c$,
  we have
  \begin{align}
    e \in Q^f(a) &\Leftrightarrow a \in \II_1;\label{ev_Q}\\
    e \in R^f(a, b) &\Leftrightarrow a \in \II_1 \text{ and } \varphi_a(b) \text{ is defined};\label{ev_R}\\
    e \in P^f(a, b, c) &\Leftrightarrow a \in \II_1, \varphi_a(b) \text{ is defined, and } \varphi_a(b)=c.\label{ev_P}
  \end{align}
  Thus each set from $Q^f(a), R^f(a, b), P^f(a, b, c)$ is $\NN$ or $\varnothing$.
  Formula \eqref{eq:contr_exam} has the form $L \to R$.
  Let us prove that $e \rvf L$ for some natural number $e$.
  It follows from (\cb), (\cek), and (\cfp) that
  there exists a $V$-function $\kk$ such that, for every $a \in \II_1$, we have $\kk(a) \in \II_2$ and
  \begin{equation}\label{L_def_k}
    \varphi_{\kk(a)}(y, y_0) \simeq \cc(\varphi_a(y), 0).
  \end{equation}
  Let $a_0 \rvf Q(a)$ for some natural numbers $a, a_0$.
  Suppose $b_0 \rvf R(a, b)$ for some natural numbers $b, b_0$.
  Using \eqref{ev_Q}, \eqref{ev_R}, we obtain $a \in \II_1$ and $!\varphi_a(b)$.
  Hence it follows from \eqref{ev_P} that $0 \rvf P(a, b, \varphi_a(b))$.
  Therefore $\cc(\varphi_a(b), 0) \rvf \exists z\,P(a, b, z)$.
  Using \eqref{L_def_k}, we get
  \begin{equation}\label{eq_Lk1}
    \varphi_{\kk(a)}(b, b_0) \rvf \exists z\,P(a, b, z).
  \end{equation}
  Thus for all natural $b, b_0$ we have \eqref{eq_Lk1} whenever $b_0 \rvf R(a, b)$.
  Therefore
  \begin{equation}\label{eq_Lk2}
    \kk(a) \rvf \forall y\,(R(a, y) \to \exists z\,P(a, y, z)).
  \end{equation}
  By (\cfp) there is a natural number $e$ such that we have $\varphi_e(x, x_0) \simeq \kk(x)$.
  It follows from \eqref{eq_Lk2} that
  \begin{equation}\label{eq_Lk3}
    \varphi_e(a, a_0) \rvf \forall y\,(R(a, y) \to \exists z\,P(a, y, z)).
  \end{equation}
  Thus for all natural $a, a_0$ we have \eqref{eq_Lk3} whenever $a_0 \rvf Q(a)$.
  Therefore $e \rvf L$.
  Using $\mathrm{(iii)}$, we get $e' \rvf R$ for some natural number $e'$.
  Thus, for all natural numbers $b, a, d_0$, we have
  $$\varphi_{e'}(b, a, d_0) \rvf \exists z\,P(a, b, z)$$ whenever $d_0 \rvf Q(a) \land R(a, b)$.
  Using \eqref{ev_Q}, \eqref{ev_R}, \eqref{ev_P} we get
  \begin{equation}\label{eq_R1}
    \pp_1 \varphi_{e'}(b, a, \cc(0, 0)) = \varphi_a(b)
  \end{equation}
  for all natural numbers $a, b$ such that $a \in \II_1$ and $!\varphi_a(b)$.
  By (\cb), (\cek), and (\cpp)
  there is a $V$-function $\uu$ such that
  \begin{equation}\label{eq_R2}
    \uu(x, y) \simeq \pp_1 \varphi_{e'}(y, x, \cc(0, 0)).
  \end{equation}
  It follows from \eqref{eq_R1}, \eqref{eq_R2} that
  \begin{equation*}
    \uu(a, b) = \varphi_a(b)
  \end{equation*}
  for all natural numbers $a, b$ such that $a \in \II_1$ and $!\varphi_a(b)$.
  Thus $\uu$ is an overuniversal $V$-function for the set of all unary $V$-functions and $\mathrm{(iii) \Rightarrow (iv)}$ is true.

  Suppose $\mathrm{(iv)}$ is true. This means that $V$ is an intuitionistic computability model.
  Using Proposition~\ref{pr:IPC_is_sound}, we get $\mathrm{(i)}$.
\end{proof}

\section{Concluding remarks}

Theorem~\ref{t_main} and \cite[Theorem~3.2]{kon_grr_bl} yield the next statement.
\begin{theorem}\label{t_main_2}
    Let $V$ be an basic computability model and
    there be no an overuniversal $V$-function for the set of all unary $V$-functions.
    Then
    \begin{itemize}
        \item Basic Predicate Calculus ($\BQC$) is sound with respect
    to the semantics of absolute $V$-realizability over all domains;
        \item Intuitionistic Predicate Calculus is not sound with respect
    to the semantics of absolute $V$-realizability over all domains.
    \end{itemize}
\end{theorem}
There are many $V$ such that they satisfy the conditions of Theorem~\ref{t_main_2}.
For example,
if $\cc$, $\pp_1$, $\pp_2$ are primitive recursive, then
the following sets of functions with some numbering satisfy the conditions:
\begin{itemize}
  \item the set of all primitive recursive functions;
  \item the set of all total recursive functions (see \cite{kon_grr_bl, kon_grr_il});
  \item the set of all arithmetical functions (see \cite{kon_ar_bl,kon_ar_pr});
  \item the set of all hyperarithmetical functions (see \cite{kon_plisko_hr}).
\end{itemize}
It is of interest to compare this models of $V$-realizability with each other.

%

\begin{acks}
This research was partially supported by Russian Foundation for Basic Research under grant 20-01-00670.
\end{acks}

\bibliographystyle{ACM-Reference-Format}
\bibliography{lit-base}


\end{document}